# The Generalized Power Generalized Weibull Distribution: Properties and Applications

Mahmoud Ali Selim

Department of Statistics, Faculty of Commerce, Al-Azher University, Egypt & King Khalid University, Community College, Saudi Arabia

**Selim.one@gmail.com**

**Abstract**

This paper introduces a new generalization of the power generalized Weibull distribution called the generalized power generalized Weibull distribution. This distribution can also be considered as a generalization of Weibull distribution. The hazard rate function of the new model has nice and flexible properties and it can take various shapes, including increasing, decreasing, upside-down bathtub and bathtub shapes. Some of the statistical properties of the new model, including quantile function, moment generating function, reliability function, hazard function and the reverse hazard function are obtained. The moments, incomplete moments, mean deviations and Bonferroni and Lorenz curves and the order statistics densities are also derived. The model parameters are estimated by the maximum likelihood method. The usefulness of the proposed model is illustrated by using two applications of real-life data.

**Keywords:** power generalized Weibull distribution, maximum likelihood estimation, moment, hazard rate function, order statistics, Bonferroni and Lorenz curves

## 1. Introduction

The three-parameter power generalized Weibull distribution was originally proposed by Bagdonavicius and Nikulin (2001) as a generalization of Weibull distribution by introducing an additional shape parameter. The cumulative density function (cdf) and probability density function (pdf) of the power generalized Weibull distribution are

$$G(x) = 1 - exp\{1 - (1 + \lambda x^\alpha)^\theta\}, \quad \alpha, \lambda, \theta > 0, \ x > 0 \quad (1)$$

and

$$g(x) = \lambda \alpha \theta x^{\alpha-1}(1 + \lambda x^\alpha)^{\theta-1} exp\{1 - (1 + \lambda x^\alpha)^\theta\}, \ \alpha, \lambda, \theta > 0, \ x > 0 \quad (2)$$

respectively, where α and θ are two shape parameters and λ is a scale parameter. The standard Weibull distribution is a special case of (1) when $\theta = 1$. This distribution can be also considered as an extension of the Nadarajah-Haghighi (NH) and exponential distributions. Nikulin and Haghighi (2007) pointed out that the hazard rate function of the power generalized Weibull distribution has nice and flexible properties and can be constant, monotone and non-monotone shaped. This distribution is often used for constructing accelerated failures times models. They also used the chi-square goodness-of-fit test to illustrate that the $PGW$ fits the randomly censored survival times data for patients at arm A of the head-and-neck cancer clinical trial. Nikulin and Haghighi (2009) have studied some of the statistical properties of PGW distribution.

Gupta et al. (1998) introduced one of the most famous and oldest methods to generalize the probability distributions called exponentiated method. If $\bar{G}(x) = [1 - G(x)]$ and $G(x)$ are the survival and cumulative density functions of the baseline distribution, then the cdf of exponentiated distribution family of Lehmann type II is defined by taking one minus the b-th -power of $\bar{G}(x)$ as follows

$$F(x) = 1 - \bar{G}(x)^b, \quad (3)$$

and the corresponding probability density function (pdf) is

$$f(x) = bg(x)\bar{G}(x)^{b-1}, \quad (4)$$

where g(x) is pdf of the baseline distribution, $b$ is a positive real parameter. Several of the generalized distributions from (3) were studied in the literature including, the generalized inverse Weibull distribution by De Gusmao et al. (2011) and the generalized inverse generalized Weibull distribution by Jain et al. (2014). The generalization of the probability distributions makes them more flexible and more suited for modeling data. Therefore, some of the generalizations of the power generalized Weibull distribution were proposed in recent years, such as; Kumaraswamy generalized power Weibull ($KGPW$) by Selim and Badr (2016) and exponentiated power generalized Weibull distribution ($EPGW$) by Pena-Ramıreza et al. (2018).

This paper aims to introduce a new generalization of the power generalized Weibull distribution named generalized power generalized Weibull ($GPGW$) distribution and studies its statistical properties. The rest of this paper is organized as follow. In section 2, the generalized power generalized Weibull distribution is introduced. In section 3, the statistical properties of the new distribution are discussed. In section 4, the reliability functions are derived. In section 5, the pdf of order statistics of the $GPGW$ model is introduced. In section 6, the maximum likelihood estimation is investigated. In section 7, two real data sets are used to illustrate the usefulness and applicability of the $GPGW$ model. In section 7, the concluding comments are given.

## 2. The Generalized Power Generalized Weibull Distribution

The reliability function $R(t)$ of the generalized power generalized Weibull distribution is simply the b-th power of the reliability function of power generalized Weibull distribution as follows

$$R(x) = \left[exp\{1 - (1 + \lambda x^\alpha)^\theta\}\right]^b, \quad \alpha, \lambda, \theta, b > 0, \ x \geq 0 \tag{5}$$

then, the cumulative density function (cdf) of the $GPGW$ distribution is

$$F(x) = 1 - R(t) = 1 - exp\{b[1 - (1 + \lambda x^\alpha)^\theta]\}, \quad \alpha, \lambda, \theta, b > 0, \ x > 0 \tag{6}$$

and the corresponding probability density function (pdf) is

$$f(x) = \alpha\lambda\theta b x^{\alpha-1}(1 + \lambda x^\alpha)^{\theta-1} exp\{b[1 - (1 + \lambda x^\alpha)^\theta]\} \quad \alpha, \lambda, \theta, b > 0, \ x > 0 \tag{7}$$

where $\lambda$ and $b$ are two scale parameters; $\alpha$ and $\theta$ are two shape parameters. The Weibull distribution is a special case of (6) when $\theta = b = 1$, hence it can be also considered as a generalization of Weibull distribution. The pdf and cdf for selected values of the parameters $\alpha$, $\theta$ and $b$ are plotted in Fig. 1, 2 respectively.

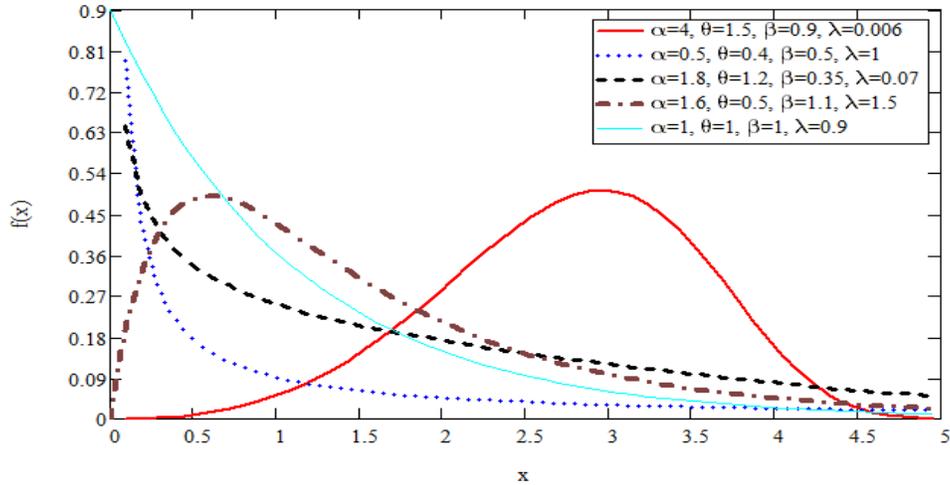

**Fig. 1: Some possible shapes of the GPGW density function**

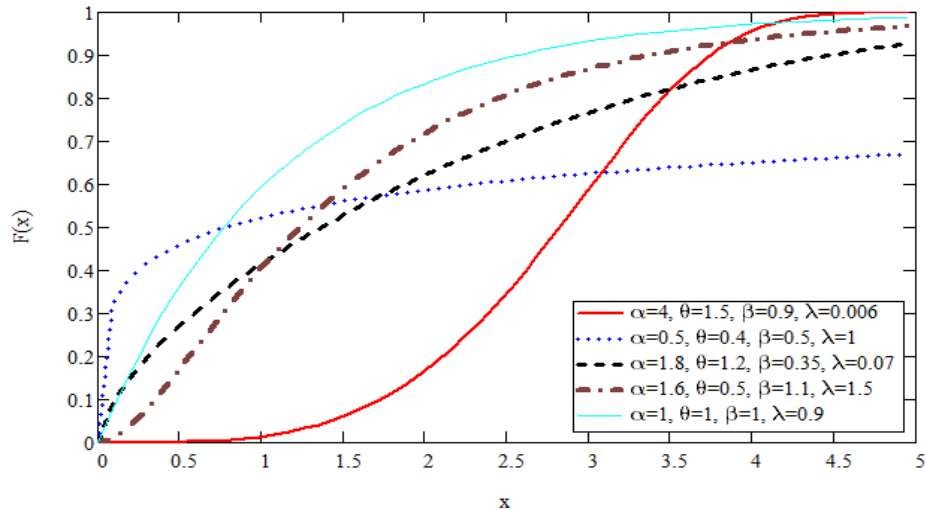

**Fig. 2: Some possible shapes of the $GPGW$ cumulative density function**

**2.1 Special distributions of the $GPGW$ distribution**

The $GPGW$ distribution having as sub-models; Weibull ($W$), exponential ($E$), power generalized Weibull distribution ($PGW$), Nadarajah-Haghighi ($NH$), Rayleigh ($R$) and generalized Nadarajah-Haghighi distribution ($GNH$) distributions. Sub-models of $GPGW$ distribution for selected values of the parameters are listed in Table 1.

**Table 1. Sub-models of the $GPGW$ distribution**

| Model | $\lambda$ | $\alpha$ | $\theta$ | $B$ | Author(s) |
|---|---|---|---|---|---|
| GPGW | - | - | - | 1 | Nikulin and Haghighi (2007) |
| W | - | - | 1 | 1 | Weibull (1951) |
| NH | - | 1 | - | 1 | Nadarajah and Haghighi (2011) |
| R | - | 2 | 1 | 1 | Rayleigh (1880) |
| E | - | 1 | 1 | 1 | Exponential distribution |
| E | 1 | 1 | 1 | - | Exponential distribution |
| GNH | - | 1 | - | - | New |

**3. The Statistical Properties of GPGW Distribution**

In this section, we derive some of the distributional properties of $GPGW$ distribution including, the quantile function, the moments, moment generating function, skewness, kurtosis and random variables generation function, incomplete moments, mean deviations and mean deviations.

**3.1 Quantile function**

The quantile function is one of the ways of specifying the distribution of a random variable and it is also an alternative to the pdf and cdf. The quantile function is usually used to obtain the statistical measures such as the median, skewness, kurtosis and it is also used to generate the random variables. The definition of the q-th quantile is the real solution of the following equation

$$F(x_q) = q, \qquad 0 \leq q \leq 1$$

Thus, the quantile function Q(q) corresponding of the $GPGW$ distribution is

$$Q(q) = \lambda^{-\frac{1}{\alpha}} \left\{ \left[ 1 - \frac{\ln(1-q)}{b} \right]^{1/\theta} - 1 \right\}^{1/\alpha} \qquad (8)$$

The median $M(x)$ of $GPGW$ distribution can be obtained from the previous function, by setting $q = 0.5$, as follows

$$M(x) = \lambda^{-\frac{1}{\alpha}}\left\{\left[1 - \frac{\ln(0.5)}{b}\right]^{1/\theta} - 1\right\}^{1/\alpha} \tag{9}$$

Fig. 3, shows the median of the $GPGW$ distribution as a function of the shape parameters $\alpha$ and $\theta$.

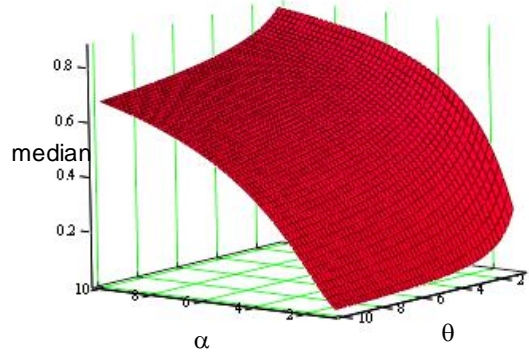

**Fig. 3: The median of the $GPGW$ distribution**

**3.2 Skewness and kurtosis**

The skewness and kurtosis measures are used in statistical analyses to characterize a distribution or a data set. The Bowley's skewness measure based on quartiles ((Kenney and Keeping 1962)) is given by

$$Sk = \frac{Q(3/4) - 2Q(1/2) + Q(1/4)}{Q(3/4) - Q(1/4)} \tag{10}$$

and the Moors' kurtosis measure based on octiles (Moors (1988)) is given by

$$Ku = \frac{Q(7/8) - Q(5/8) + Q(3/8) - Q(1/8)}{Q(6/8) - Q(2/8)} \tag{11}$$

The skewness and kurtosis measures based on quantiles like Bowley's skewness and Moors' kurtosis have a number of advantages compared to the classical measures of skewness and kurtosis, e.g. they are less sensitive to outliers and they exist for the distributions even without defined the moments. Fig. 4, shows the behaviors of skewness and kurtosis of the $GPGW$ distribution as a function of the shape parameters $\alpha$ and $\theta$.

**3.3 Random variables generation**

The closed form of the quantile function of the $GPGW$ distribution makes the simulation from this distribution easier. Therefore, the random variables of $GPGW$ distribution are directly generated from the following function

$$X = \lambda^{-\frac{1}{\alpha}}\left\{\left[1 - \frac{\ln(1-u)}{b}\right]^{1/\theta} - 1\right\}^{1/\alpha} \tag{12}$$

where $\alpha$, $\lambda$, $\theta$ and $b$ are known parameters and $u$ is generated number from the uniform distribution $(0, 1)$. Now, we use eq (12) to generate two simulated data sets for some parameter values. Fig. 5, displays the histograms and the exact $GPGW$ densities for the simulated data sets. These shapes show that the setting is quite adequate and reinforces that this model has good potential for simulation studies.

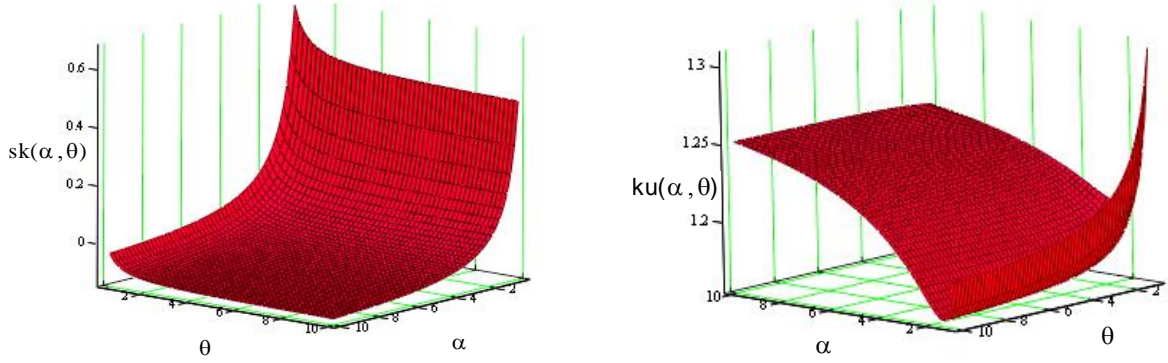

**Fig. 4:** The Skewness (sk) and kurtosis (ku) of the $GPGW$ distribution

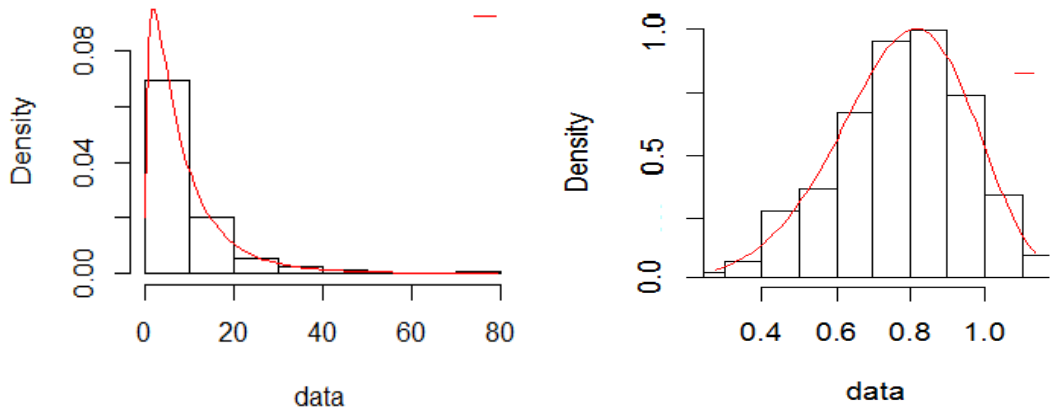

**Fig. 5:** Plots of the $GPGW$ pdf and histogram for simulated data with n = 300 and n=2000 respectively

### 3.4 The Moments

If $X$ has the $GPGW$ distribution, then the r-th moment of $X$ for integer value of $r/\alpha$ is

$$\mu'_r = be^b \lambda^{-\frac{r}{\alpha}} \sum_{j=0}^{r/\alpha} \frac{(-1)^{j+\frac{r}{\alpha}}}{b^{\frac{j}{\theta}+1}} \binom{r/\alpha}{j} \Gamma\left(\frac{j}{\theta}+1, b\right) \quad (13)$$

**Proof.** The r-th moment is defined as follows

$$\mu'_r = E(X^r) = \int_0^\infty x^r f(x) dx \quad (14)$$

substituting $f(x)$ from (7) in (14) yields

$$\mu'_r = \alpha\theta\lambda b \int_0^\infty x^{r+\alpha-1}(1+\lambda x^\alpha)^{\theta-1} e^{b(1-(1+\lambda x^\alpha)^\theta)} dx, \quad (15)$$

by taking $u^{\frac{1}{\theta}} = 1 + \lambda x^\alpha$, the above expression reduce to

$$\mu'_r = be^b \lambda^{-\frac{r}{\alpha}} \int_1^\infty \left(u^{\frac{1}{\theta}} - 1\right)^{\frac{r}{\alpha}} e^{-bu} du, \quad (16)$$

by applying the binomial expansion, then (16) become

$$\mu_r' = b e^b \lambda^{-\frac{r}{\alpha}} \sum_{j=0}^{r/\alpha} (-1)^{\frac{r}{\alpha}+j} \binom{r/\alpha}{j} \int_1^\infty u^{\frac{j}{\theta}} e^{-bu} du, \tag{17}$$

by integrating the incomplete gamma function $\int_1^\infty u^{\frac{j}{\theta}} e^{-bu} du = \Gamma\left(\frac{j}{\theta}+1, b\right) / b^{\frac{j}{\theta}+1}$, we get the r-*th* moment of X as follows

$$\mu_r' = b e^b \lambda^{-\frac{r}{\alpha}} \sum_{j=0}^{r/\alpha} \frac{(-1)^{j+\frac{r}{\alpha}}}{b^{\frac{j}{\theta}+1}} \binom{r/\alpha}{j} \Gamma\left(\frac{j}{\theta}+1, b\right). \quad \blacksquare$$

If $b=1$, we get the moments of Nikulin- Haghighi distribution as follows

$$\mu_r' = \lambda^{-\frac{r}{\alpha}} e \sum_{j=0}^{r/\alpha} (-1)^{\frac{r}{\alpha}+j} \binom{r/\alpha}{j} \Gamma\left(\frac{j}{\theta}+1, 1\right)$$

which agrees with Nikulin and Haghighi (2007). Also, if $\alpha = b = 1$, we get the moments of Nadarajah-Haghighi distribution as follows

$$\mu_r' = \lambda^{-r} e \sum_{j=0}^{r} (-1)^{r+j} \binom{r}{j} \Gamma\left(\frac{j}{\theta}+1, 1\right)$$

which agrees with the moments obtained by Nadarajah and Haghighi (2011).

In particular, the first and second moments and the variance of X can be determined from (13) as follow

$$\mu_1' = E(X) = b e^b \lambda^{-\frac{1}{\alpha}} \sum_{j=0}^{1/\alpha} \frac{(-1)^{j+\frac{1}{\alpha}}}{b^{\frac{j}{\theta}+1}} \binom{1/\alpha}{j} \Gamma\left(\frac{j}{\theta}+1, b\right), \tag{18}$$

$$\mu_2' = E(X^2) = b e^b \lambda^{-\frac{2}{\alpha}} \sum_{j=0}^{2/\alpha} \frac{(-1)^{j+\frac{2}{\alpha}}}{b^{\frac{j}{\theta}+1}} \binom{2/\alpha}{j} \Gamma\left(\frac{j}{\theta}+1, b\right) \tag{19}$$

and

$$Var(X) = \mu_2' - [\mu_1']^2 \tag{20}$$

respectively. The non- central moments in (13) can be also used to calculated the central moments $\mu_r$ and the cumulants $k_r$ as following $\mu_r = \sum_{k=0}^{r}(-1)^k \binom{r}{k} \mu_1'^k \mu_{r-k}'$ and $k_r = \mu_r' - \sum_{k=1}^{r-1} \binom{r-1}{k-1} k_r \mu_{r-k}'$, where the cumulants $k_r$ are defined as quantities that provide an alternative to the distribution moments. The skewness $\gamma_1$ and kurtosis $\gamma_2$ can be calculated based on cumulants in the forms $\gamma_1 = k_3/k_2^{3/2}$ and $\gamma_2 = k_4/k_2^2$, respectively.

### 3.5 The moment generating function

If X~GPGW distribution, then for any integer value of $r/\alpha$ the moment generating function is

$$M_x(t) = b e^b \lambda^{-\frac{r}{\alpha}} \sum_{r=0}^{\infty} \sum_{j=0}^{\infty} \frac{(-1)^{j+\frac{r}{\alpha}} t^r}{b^{\frac{j}{\theta}+1} r!} \binom{r/\alpha}{j} \Gamma\left(\frac{j}{\theta}+1, b\right) \tag{21}$$

**Proof.** The moment generating function is defined as follows

$$M_x(t) = \int_0^\infty e^{tx} f(x) dx$$

Using exponential function formula $e^{tx} = \sum_{r=0}^{\infty} \frac{(tx)^r}{r!}$, we get

$$M_x(t) = \sum_{r=0}^{\infty} \frac{t^r}{r!} \int_0^\infty x^r f(x) dx = \sum_{r=0}^{\infty} \frac{t^r}{r!} \mu_r' \tag{22}$$

Substituting (16) in to (22) yields the mgf of $GPGW$ distribution as in (21). $\blacksquare$

### 3.6 Incomplete moments

The sth incomplete moment of $X$, is defined as follows

$$m_r(s) = E(X^r | X > s) = \int_s^\infty x^r f(x) dx \tag{23}$$

Hence, by inserting (7) in (23) and after some manipulate, we get the sth incomplete moment of $GPGW$ distribution as follows

$$m_r(s) = b e^b \lambda^{-\frac{r}{\alpha}} \sum_{j=0}^{r/\alpha} \frac{(-1)^{j+\frac{r}{\alpha}}}{b^{\frac{j}{\theta}+1}} \binom{r/\alpha}{j} \Gamma\left(\frac{j}{\theta} + 1, b(1 + \lambda s^\alpha)^\theta\right) \tag{24}$$

**Proof.** Substituting (7) into (23), yields

$$m_r(s) = \alpha \theta \lambda b \int_s^\infty x^{r+\alpha-1} (1 + \lambda x^\alpha)^{\theta-1} e^{b(1-(1+\lambda x^\alpha)^\theta)} dx \tag{25}$$

By setting $u = (1 + \lambda x^\alpha)^\theta$, the above expression reduce to

$$m_r(s) = b e^b \lambda^{-\frac{r}{\alpha}} \int_{(1+\lambda s^\alpha)^\theta}^\infty \left(u^{\frac{1}{\theta}} - 1\right)^{\frac{r}{\alpha}} e^{-bu} du \tag{26}$$

Applying the binomial expansion gives

$$m_r(s) = b \lambda^{-\frac{r}{\alpha}} e^b \sum_{j=0}^{r/\alpha} (-1)^{\frac{r}{\alpha}+j} \binom{r/\alpha}{j} \int_{(1+\lambda s^\alpha)^\theta}^\infty u^{\frac{j}{\theta}} e^{-bu} du \tag{27}$$

By integrating the incomplete gamma function $\int_{(1+\lambda s^\alpha)^\theta}^\infty u^{\frac{j}{\theta}} e^{-bu} du = \Gamma\left(\frac{j}{\theta} + 1, b(1 + \lambda s^\alpha)^\theta\right)/b^{\frac{j}{\theta}+1}$, we get the $r$ th upper incomplete moment of $GPGW$ distribution as follows

$$m_r(s) = b \lambda^{-\frac{r}{\alpha}} e^b \sum_{j=0}^{r/\alpha} \frac{(-1)^{j+\frac{r}{\alpha}}}{b^{\frac{j}{\theta}+1}} \binom{r/\alpha}{j} \Gamma\left(\frac{j}{\theta} + 1, b(1 + \lambda s^\alpha)^\theta\right). \blacksquare$$

In particular, the first incomplete moments of the $GPGW$ distribution can be obtained by putting $r = 1$ in (24), as follows

$$m_1(s) = b e^b \lambda^{-\frac{1}{\alpha}} \sum_{j=0}^{1/\alpha} \frac{(-1)^{j+\frac{1}{\alpha}}}{b^{\frac{j}{\theta}+1}} \binom{1/\alpha}{j} \Gamma\left(\frac{j}{\theta} + 1, b(1 + \lambda s^\alpha)^\theta\right) \tag{28}$$

### 3.7 Mean deviations

The mean deviations about the mean ($\delta_1(X)$) and about the median ($\delta_1(X)$) of $X$ are given, respectively, by

$$\delta_1(X) = E(|X - \mu_1'|) = \int_0^\infty |X - \mu_1'| f(x) dx$$

$$= 2\mu_1' F(\mu_1') - 2\mu_1' + 2 \int_{\mu_1'}^\infty x f(x) dx$$

$$= 2\mu_1' F(\mu_1') - 2H_1(\mu_1') \tag{29}$$

and

$$\delta_2(X) = E(|X - M|) = \int_0^\infty |X - M| f(x) dx$$

$$= 2MF(M) - M - \mu_1' + 2 \int_M^\infty x f(x) dx$$

$$= \mu_1' - 2H_1(M) \tag{30}$$

where $\mu_1' = E(X)$, $M$ = median(X), $F(\mu_1')$ from (6) and $H_1(s)$ is the sth lower incomplete moment as follows

$$H_1(s) = \int_0^{\mu_1'} x f(x) dx$$

$$= b e^b \lambda^{-\frac{1}{\alpha}} \sum_{j=0}^{1/\alpha} \frac{(-1)^{j+\frac{1}{\alpha}}}{b^{\frac{j}{\theta}+1}} \binom{r/\alpha}{j} \left[ \Gamma\left(\frac{j}{\theta}+1, b\right) - \Gamma\left(\frac{j}{\theta}+1, b(1+\lambda s^\alpha)^\theta\right) \right] \quad (31)$$

**3.8 Bonferroni and Lorenz curves**

The Bonferroni and Lorenz curves (Bonferroni (1930)) have many applications, especially in economics to study income and poverty and in other fields like reliability, demography, insurance and medicine. The Bonferroni and Lorenz curves are defined as

$$B(p) = \frac{1}{p\mu_1'} \int_0^q x f(x) dx = \frac{\mu_1' - m_1(q)}{p\mu_1'} \quad (32)$$

and

$$L(p) = \frac{1}{\mu_1'} \int_0^q x f(x) dx = \frac{\mu_1' - m_1(q)}{\mu_1'} \quad (33)$$

respectively, where $\mu_1' = E(X)$ and $q = Q(p)$ is calculated from (8) for a given probability $(p)$, and $m_1(q)$ is the first incomplete moment from (28).

$$B(p) = \frac{1}{p} - \frac{b e^b \lambda^{-\frac{1}{\alpha}} \sum_{j=0}^{1/\alpha} \frac{(-1)^{j+\frac{1}{\alpha}}}{b^{\frac{j}{\theta}+1}} \binom{1/\alpha}{j} \left[ \Gamma\left(\frac{j}{\theta}+1, b(1+\lambda q^\alpha)^\theta\right) \right]}{p\mu_1'} \quad (34)$$

and

$$L(p) = 1 - \frac{b e^b \lambda^{-\frac{1}{\alpha}} \sum_{j=0}^{1/\alpha} \frac{(-1)^{j+\frac{1}{\alpha}}}{b^{\frac{j}{\theta}+1}} \binom{1/\alpha}{j} \left[ \Gamma\left(\frac{j}{\theta}+1, b(1+\lambda q^\alpha)^\theta\right) \right]}{\mu_1'} \quad (35)$$

**4. Reliability Analysis**

In this section, the survival $s(t)$, failure rate $h(t)$, reversed hazard $r(t)$ and the cumulative failure rate $H(t)$ functions of $GPGW$ distribution are derived.

**4.1 The survival function**

The survival function $R(t)$ of the $GPGW$ distribution can be derived using the cumulative distribution function in (5) as follows

$$R(t) = 1 - F(x) = \left[ e^{1-(1+\lambda x^\alpha)^\theta} \right]^b, \quad t > 0 \quad (36)$$

**4.2 The hazard function**

For a continuous distribution with pdf $f(x)$ and cdf $F(x)$, the hazard rate function for any time is defined as follows

$$h(t) = \lim_{\Delta t \to 0} \frac{P(T < t + \Delta x | T > t)}{\Delta t} = \frac{f(t)}{1 - F(t)}$$

Subsequently, the hazard rate for any time of the $GPGW$ distribution can be determined using the cdf and pdf in Eqs. (5), (6) as follow

$$h(t) = \alpha \lambda \theta b t^{\alpha-1} (1 + \lambda t^\alpha)^{\theta-1}, \quad t > 0 \quad (37)$$

Nikulin and Haghighi (2007) discussed and established that the hazard rate function of the $PGW$ distribution is increasing, decreasing, bathtub and unimodal failure rate. It follows that the hazard rate function for $GPGW$ distribution can be taking the following forms:

- monotone increasing if either $\alpha > 1$ and $\alpha\theta > 1$ or $\alpha = 1$ and $\theta > 1$;
- monotone decreasing if either $0 < \alpha < 1$ and $\alpha\theta < 1$ or $0 < \alpha < 1$ and $\alpha\theta = 1$;
- bathtub shaped if $0 < \alpha < 1$ and $\alpha\theta > 1$;
- unimodal (inverted bathtub shaped) if $\alpha > 1$ and $0 < \alpha\theta < 1$;
- constant, $h(t) = b\lambda$ if $\alpha = \theta = 1$.

The plots of the hazard function of $GPGW$ distribution for some selected parameters values are displayed in Fig. 6. These plots show the flexibility of hazard rate function that makes the $GPGW$ hazard rate function useful and suitable for non-monotone hazard behaviors that are more likely to be observed in real life situations.

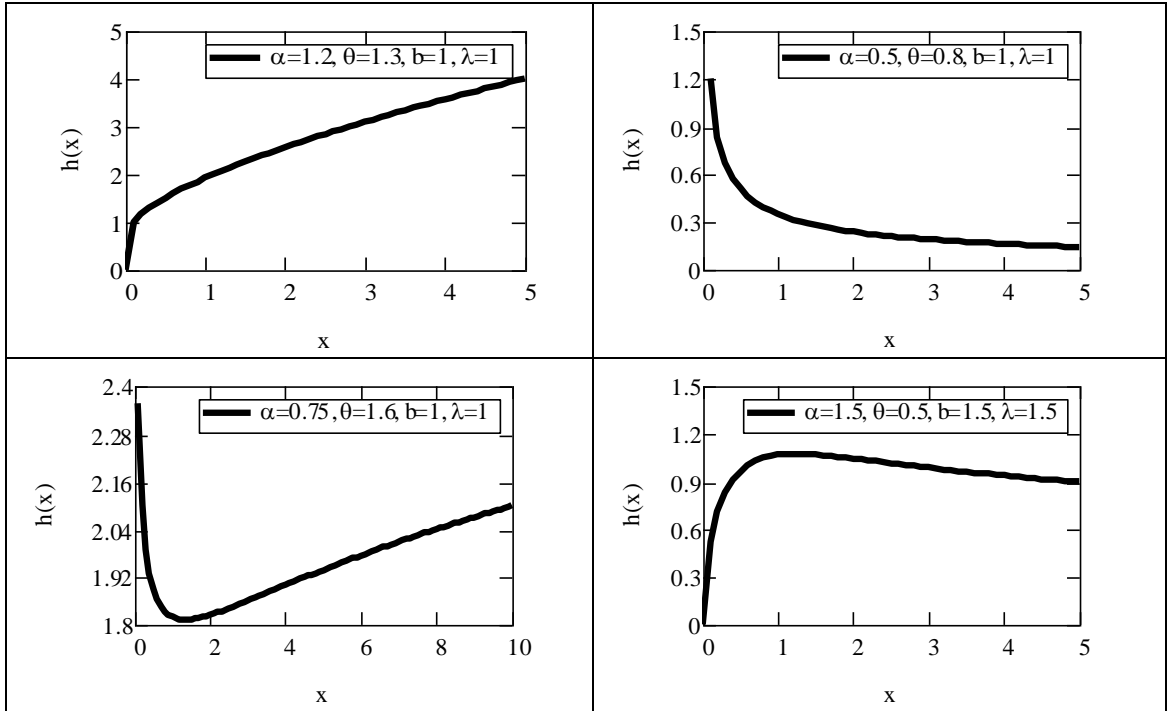

**Fig. 6: Some possible shapes of the $GPGW$ hazard rate function**

**4.3 The reversed hazard and cumulative hazard rate functions**

The reversed hazard $r(t)$ and the cumulative hazard rate $H(t)$ functions of $GPGW$ distribution are given, respectively, as follow

$$r(t) = \frac{\alpha\lambda\theta b x^{\alpha-1}(1+\lambda x^\alpha)^{\theta-1}\left[e^{1-(1+\lambda x^\alpha)^\theta}\right]^b}{1 - \left[e^{1-(1+\lambda x^\alpha)^\theta}\right]^b}, \quad t > 0, \tag{38}$$

and

$$H(t) = -\ln\left[1 - \exp\{b[1 - (1+\lambda x^\alpha)^\theta]\}\right], \quad t > 0. \tag{39}$$

**5. Order Statistics**

The order statistics arise naturally in many areas of statistical theory and practice which makes it one of the important statistical topics. Let $X_{(1)}, X_{(2)}, \dots, X_{(n)}$ denote the order statistics of a random sample drawn from a continuous distribution with cdf $F(x)$ and pdf $f(x)$, then the pdf of $X_{(k)}$ is given by

$$f_{k:n}(x) = \frac{n!}{(k-1)!\,(n-k)!} f(x)[F(x)]^{k-1}[1-F(x)]^{n-k}, \quad k = 1,2,\dots,n \tag{40}$$

Let $X$ is a random variable of $GPGW$ distribution, then by substituting (6) and (7) into equation (40), we get the $kth$ order statistics of $GPGW$ density function as follows

$$f_{k:n}(x) = \frac{n!}{(k-1)!\,(n-k)!} \alpha\lambda\theta b x^{\alpha-1}(1+\lambda x^\alpha)^{\theta-1}\left[1-\left(e^{1-(1+\lambda x^\alpha)^\theta}\right)^b\right]^{k-1}$$
$$\times \left[e^{1-(1+\lambda x^\alpha)^\theta}\right]^{b(n-k+1)} \tag{41}$$

The pdf of order statistics when $k = 1$ and when $k = n$ are

$$f_{1:n}(x) = n\alpha\lambda\theta b x^{\alpha-1}(1+\lambda x^\alpha)^{\theta-1}\left[e^{1-(1+\lambda x^\alpha)^\theta}\right]^{nb} \tag{42}$$

and

$$f_{n:n}(x) = n\alpha\lambda\theta b x^{\alpha-1}(1+\lambda x^\alpha)^{\theta-1}\left[1-\left(e^{1-(1+\lambda x^\alpha)^\theta}\right)^b\right]^{n-1}\left[e^{1-(1+\lambda x^\alpha)^\theta}\right]^b \tag{43}$$

respectively.

## 6. Maximum Likelihood Estimation

This section discusses the maximum likelihood estimation ($MLE$) for the parameters $\sigma = (\alpha, \lambda, \theta, b)$ of the $GPGW$ distribution. Let $x_1, x_1, \ldots, x_n$ is a complete random sample of size $n$ from the $GPGW$ distribution. Then the likelihood function (LF) is

$$L(\sigma|x) = (\alpha\lambda\theta b)^n \prod_{i=1}^{n} x_i^{\alpha-1}(\lambda x_i^\alpha + 1)^{\theta-1}\left(e^{1-(\lambda x_i^\alpha+1)^\theta}\right)^b \tag{44}$$

and the log-likelihood function ($logL$) is

$$logL = n\log(\alpha\lambda\theta b) + (\alpha-1)\sum_{i=1}^{n}\ln x_i + (\theta-1)\sum_{i=1}^{n}\ln(\lambda x_i^\alpha + 1) + b$$
$$- b\sum_{i=1}^{n}(\lambda x_i^\alpha + 1)^\theta \tag{45}$$

The partial derivatives of the above equation are

$$\frac{\partial \ln L}{\partial b} = \frac{n}{b} - \sum_{i=1}^{n}(\lambda x_i^\alpha + 1)^\theta + 1, \tag{46}$$

$$\frac{\partial \ln L}{\partial \alpha} = \frac{n}{\alpha} + \sum_{i=1}^{n}\ln(x_i) - \theta\lambda b\sum_{i=1}^{n}\ln(x_i)x_i^\alpha(\lambda x_i^\alpha + 1)^{\theta-1} + (\theta-1)\lambda\sum_{i=1}^{n}\frac{\ln(x_i)x_i^\alpha}{\lambda x_i^\alpha + 1}, \tag{47}$$

$$\frac{\partial \ln L}{\partial \lambda} = \frac{n}{\lambda} - \theta b\sum_{i=1}^{n}x_i^\alpha(\lambda x_i^\alpha + 1)^{\theta-1} + (\theta-1)\sum_{i=1}^{n}\frac{x_i^\alpha}{\lambda x_i^\alpha + 1}, \tag{48}$$

$$\frac{\partial \ln L}{\partial \theta} = \frac{n}{\theta} + \sum_{i=1}^{n}\ln(\lambda x_i^\alpha + 1) - b\sum_{i=1}^{n}\ln(\lambda x_i^\alpha + 1)(\lambda x_i^\alpha + 1)^\theta. \tag{49}$$

The maximum likelihood estimators of $\alpha, \lambda, \theta$ and $b$ are the simultaneous solutions of the following nonlinear likelihood equations

$$\frac{\partial \ln L}{\partial b} = \frac{\partial \ln L}{\partial \alpha} = \frac{\partial \ln L}{\partial \lambda} = \frac{\partial \ln L}{\partial \theta} = 0$$

The previous equations cannot be analytically solved, but the statistical software can be used to solve them numerically by using iterative techniques like the Newton-Raphson algorithm.

**7. Real Data Illustration**

In this section, the usefulness of the $GPGW$ distribution is illustrated by using two real datasets. These datasets are described as follows

**The dataset (I): Failure times of 50 devices**

The first data set represents the failure times of 50 devices put under a life test (see Aarset (1987)). These data are: 0.1, 0.2, 1.0, 1.0, 1.0, 1.0, 1.0, 2.0, 3.0, 6.0, 7.0, 11.0, 12.0, 18.0, 18.0, 18.0, 18.0, 18.0, 21.0, 32.0, 36.0, 40.0, 45.0, 45.0, 47.0, 50.0, 55.0, 60.0, 63.0, 63.0, 67.0, 67.0, 67.0, 67.0, 72.0, 75.0, 79.0, 82.0, 82.0, 83.0, 84.0, 84.0, 84.0, 85.0, 85.0, 85.0, 85.0, 85.0, 86.0, 86.0 .

**The dataset (II): Survival times of acute myelogenous leukaemia data**

The second real data set represents the survival times, in weeks, of 33 patients suffering from acute myelogenous leukaemia (see Feigl and Zelen (1965)). The data are: 65, 156, 100, 134, 16, 108, 121, 4, 39, 143, 56, 26, 22, 1, 1, 5, 65, 56, 65, 17, 7, 16, 22, 3, 4, 2, 3, 8, 4, 3, 30, 4, 43.

knowing the shape of hazard function curve can help to guess the appropriate distribution of data. For this purpose, a useful device called the total time on test (TTT) plot is used to identifying the shape of hazard rate function graphically (see for more details Aarset (1987)). The TTT plot is obtained by plotting $G(r/n) = \left[\sum_{i=1}^{n} x_{(i)} + (n-r)x_{(r)}\right]/\sum_{i=1}^{n} x_{(i)}$, against $r/n$, where $i, r = 1, \ldots, n$ and $x_{(i)}$ are the order statistics of the sample. Fig. 7 displays the $TTT$ plots for our datasets. The TTT plot for dataset I in Fig. 7 (a) indicates a bathtub shaped hazard rate function. Also, the $TTT$ plot for dataset II in Fig. 7 (b) points out a decreasing hazard rate function. These shapes reveal the suitability of the $GPGW$ distribution to fit our datasets.

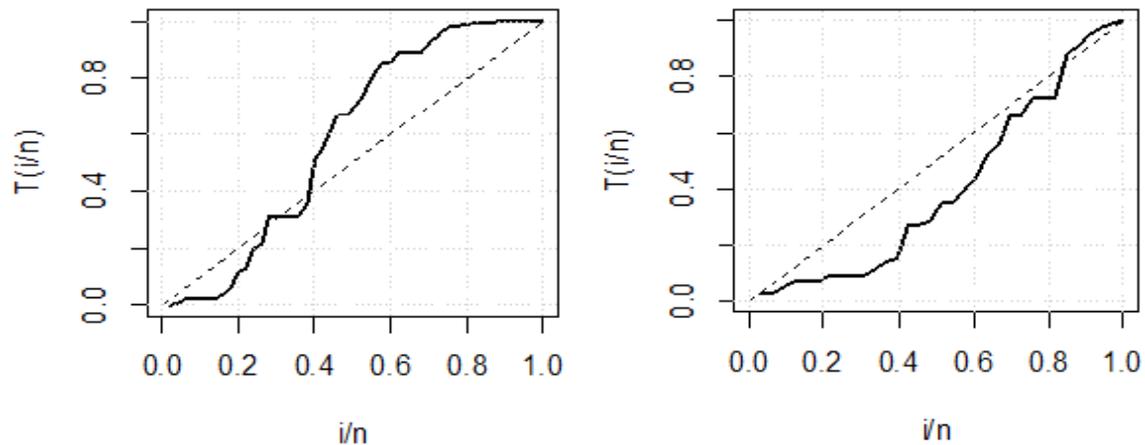

**Fig.7: (a) TTT plot for dataset I (in the left ); (b) TTT plot for dataset II (in the right)**

Now, we fit the $GPGW$ distribution to our datasets and compare it with the other fitted models such power generalized Weibull ($PGW$), exponentiated power generalized Weibull ($EPGW$), Nadarajah-Haghighi ($NH$), exponentiated Nadarajah-Haghighi ($ENH$), exponentiated Weibull ($EW$), Weibull ($W$), exponentiated exponential ($EE$) and exponential ($E$) distributions. The required numerical evaluations are implemented using the R software (Team (2013)). In order to compare the previous distributions and to verify the quality of the fits, we consider some of well-known Goodness-of-Fit statistics like, Cramér-von Mises ($W^*$), Anderson Darling ($A^*$), Kolmogorov-Smirnov ($KS$), Maximized Loglikelihood ($-L$), Akaike Information Criterion ($AIC$) and the Consistent Akaike Information Criterion ($CAIC$). The model with a minimum value of Goodness-of-Fit statistics is the best model to fit the data.

The results of the maximum likelihood estimates and their standard errors of the fitted models for dataset I and dataset II are displayed in Tables 2 and 3 respectively. As well, the values of Goodness-of-Fit statistics of the fitted models for dataset I and dataset II are displayed in Tables 4 and 5 respectively. Tables 4 and 5, show that the proposed $GPGW$ model gives smallest values for the Goodness-of-Fit statistics. The plots of histogram and estimated pdf for the $GPGW$ distribution and other fitted distributions for datasets I and II are displayed in Fig. 8 and Fig. 9, respectively. These plots also indicate that the $GPGW$ distribution provides an adequate fit than other distributions for both data sets.

**Table 2. The estimates and the standard errors (in parentheses) for dataset I**

| Model | Estimates | | | |
|---|---|---|---|---|
| | $\theta$ | $\lambda$ | $\alpha$ | $\beta$ |
| $GPGW$ | 10.8630 (16.2942) | 0.0232 (0.0770) | 0.5872 (0.2732) | 0.1101 (0.2031) |
| $EPGW$ | 3.4728 (2.0588) | 0.0024 (0.0006) | 1.1088 (0.1490) | 0.6557 (0.1303) |
| $PGW$ | 7.2704 (4.004) | 0.0045 (0.0022) | 0.8074 (0.0986) | - |
| $EW$ | - | 0.0029 (0.0010) | 1.3965 (0.0832) | 0.5356 (0.0876) |
| $ENH$ | 5.0595 (1.2794) | 0.0026 (0.0007) | - | 0.7243 (0.1146) |
| $NH$ | 4.2584 (1.4363) | 0.0036 (0.0013) | - | - |
| $EE$ | - | 0.0187 (0.0036) | - | 0.7801 (0.1351) |
| $W$ | - | 0.0272 (0.0137) | 0.9476 (0.1174) | - |
| $E$ | - | 0.0219 (0.0031) | - | - |

**Table 3. The estimates and the standard errors (in parentheses) for dataset II**

| Model | Estimates | | | |
|---|---|---|---|---|
| | $\theta$ | $\lambda$ | $\alpha$ | $\beta$ |
| $GPGW$ | 0.0373 (0.1342) | 7.8281 (21.8268) | 6.5271 (57.6467) | 0.1239 (0.1026) |
| $EPGW$ | 0.0385 (0.0879) | 19.5120 (65.3883) | 6.7869 (15.6964) | 2.0399 (1.2717) |
| $PGW$ | 1.4058 (3.0889) | 0.0475 (0.0905) | 0.7200 (0.3198) | - |
| $EW$ | - | 0.1680 (0.5769) | 0.5998 (0.6057) | 1.5748 (2.9553) |
| $ENH$ | 0.3151 (0.2651) | 0.6516 (3.214) | - | 1.8757 (3.5269) |
| $NH$ | 0.4897 (0.1482) | 0.0998 (0.0701) | - | - |
| $EE$ | - | 0.6785 (0.1448) | - | 0.0188 (0.0048) |
| $W$ | - | 0.0628 (0.0298) | 0.7763 (0.1073) | - |
| $E$ | - | 0.0245 (0.0043) | - | - |

Table 4. Goodness-of-fit statistics for dataset I

| Model | Kolmogorov-Smirnov | | $W^*$ | $A^*$ | -L | AIC | CAIC |
|---|---|---|---|---|---|---|---|
| | K–S | P-value | | | | | |
| GPGW | 0.1675 | 0.1209 | 0.3201 | 2.0594 | 232.732 | 473.464 | 474.353 |
| EPGW | 0.2084 | 0.0261 | 0.3716 | 2.3410 | 234.497 | 476.995 | 477.884 |
| PGW | 0.1977 | 0.0402 | 0.3952 | 2.4683 | 235.879 | 477.757 | 478.279 |
| EW | 0.2084 | 0.0260 | 0.4261 | 2.6343 | 237.311 | 480.622 | 481.144 |
| ENH | 0.2033 | 0.0320 | 0.3753 | 2.3612 | 234.5782 | 475.156 | 475.678 |
| NH | 0.1881 | 0.0581 | 0.3749 | 2.3584 | 237.182 | 478.365 | 478.620 |
| EE | 0.2043 | 0.0308 | 0.4837 | 2.9421 | 239.973 | 483.947 | 484.202 |
| W | 0.1930 | 0.0483 | 0.4948 | 3.0009 | 240.980 | 485.959 | 486.215 |
| E | 0.1914 | 0.0514 | 0.4861 | 2.9546 | 241.068 | 484.135 | 484.219 |

Table 5. Goodness-of-fit statistics for dataset II

| Model | Kolmogorov-Smirnov | | $W^*$ | $A^*$ | -L | AIC | CAIC |
|---|---|---|---|---|---|---|---|
| | K–S | P-value | | | | | |
| GPGW | 0.1289 | 0.6424 | 0.0688 | 0.4710 | 151.199 | 310.399 | 311.828 |
| EPGW | 0.1354 | 0.5805 | 0.1095 | 0.7058 | 153.466 | 314.932 | 316.361 |
| PGW | 0.1366 | 0.5694 | 0.0950 | 0.6565 | 153.571 | 313.143 | 313.970 |
| EW | 0.1366 | 0.5692 | 0.0954 | 0.6455 | 153.562 | 313.123 | 313.951 |
| ENH | 0.1301 | 0.6315 | 0.1103 | 0.6949 | 153.746 | 313.492 | 314.320 |
| NH | 0.1393 | 0.5441 | 0.1002 | 0.6659 | 153.743 | 311.486 | 311.886 |
| EE | 0.1384 | 0.5521 | 0.0966 | 0.6691 | 153.652 | 311.303 | 311.703 |
| W | 0.1366 | 0.5689 | 0.0948 | 0.6508 | 153.587 | 311.173 | 311.574 |
| E | 0.2182 | 0.0864 | 0.0973 | 0.6730 | 155.450 | 312.900 | 313.029 |

**7. Conclusion**

This paper introduces a new four-parameter distribution called the generalized power generalized Weibull distribution. This distribution has a number of well-known distributions as sub-models like Weibull, Rayleigh, exponential, power generalized Weibull, Nadarajah-Haghighi distributions. The statistical properties of the new model including the hazard rate function, quantile function, order statistics, moments, incomplete moments, mean deviations and Bonferroni and Lorenz curves are derived. The hazard function of the new model exhibits a variety of shapes, like increasing, decreasing, bathtub, unimodal and constant shapes. The maximum likelihood method is used to estimate the model parameters. The real applications have established that the new distribution is quite useful for dealing with lifetime data and behaves better than other distributions that are commonly used for fitting this type of data.

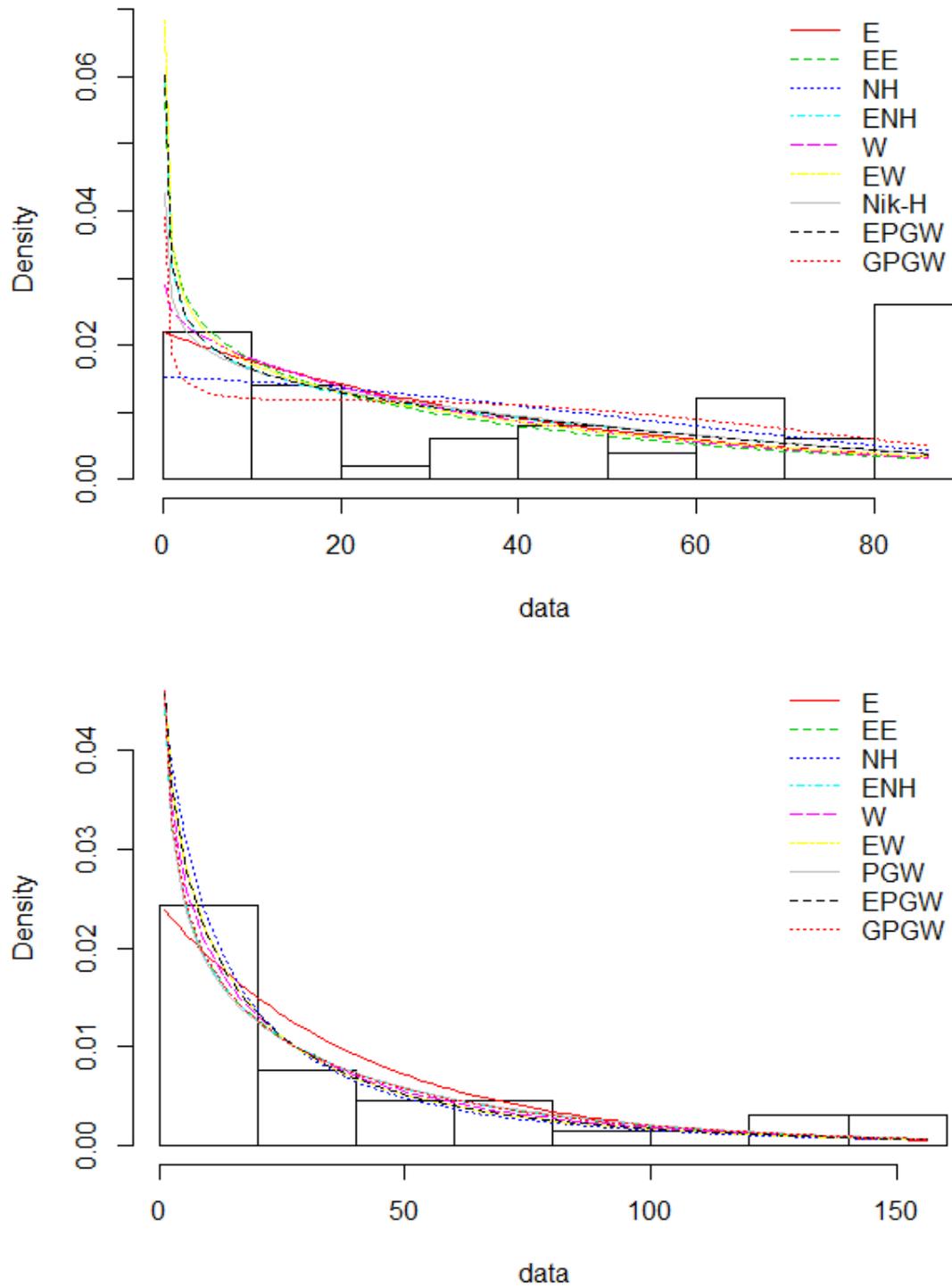

**Fig. 8: Histogram and estimated densities for dataset I (upper Panel); Histogram and estimated densities for dataset II (bottom Panel).**